\documentclass{amsart}
\usepackage{latexsym}
\usepackage{amssymb}
\usepackage{amsmath}

\begin{document}


\newtheorem{theorem}{Theorem}
\newtheorem{problem}{Problem}
\newtheorem{definition}{Definition}
\newtheorem{lemma}{Lemma}
\newtheorem{proposition}{Proposition}
\newtheorem{corollary}{Corollary}
\newtheorem{example}{Example}
\newtheorem{conjecture}{Conjecture}
\newtheorem{algorithm}{Algorithm}
\newtheorem{exercise}{Exercise}
\newtheorem{xample}{Example}
\newtheorem{remarkk}{Remark}

\newcommand{\be}{\begin{equation}}
\newcommand{\ee}{\end{equation}}
\newcommand{\bea}{\begin{eqnarray}}
\newcommand{\eea}{\end{eqnarray}}
\newcommand{\beq}[1]{\begin{equation}\label{#1}}
\newcommand{\eeq}{\end{equation}}
\newcommand{\beqn}[1]{\begin{eqnarray}\label{#1}}
\newcommand{\eeqn}{\end{eqnarray}}
\newcommand{\beaa}{\begin{eqnarray*}}
\newcommand{\eeaa}{\end{eqnarray*}}
\newcommand{\req}[1]{(\ref{#1})}

\newcommand{\lip}{\langle}
\newcommand{\rip}{\rangle}
\newcommand{\uu}{\underline}
\newcommand{\oo}{\overline}
\newcommand{\La}{\Lambda}
\newcommand{\la}{\lambda}
\newcommand{\eps}{\varepsilon}
\newcommand{\om}{\omega}
\newcommand{\Om}{\Omega}
\newcommand{\ga}{\gamma}
\newcommand{\ka}{\kappa}
\newcommand{\rrr}{{\Bigr)}}
\newcommand{\qqq}{{\Bigl\|}}

\newcommand{\dint}{\displaystyle\int}
\newcommand{\dsum}{\displaystyle\sum}
\newcommand{\dfr}{\displaystyle\frac}
\newcommand{\bige}{\mbox{\Large\it e}}
\newcommand{\integers}{{\Bbb Z}}
\newcommand{\rationals}{{\Bbb Q}}
\newcommand{\reals}{{\rm I\!R}}
\newcommand{\realsd}{\reals^d}
\newcommand{\realsn}{\reals^n}
\newcommand{\NN}{{\rm I\!N}}
\newcommand{\DD}{{\rm I\!D}}
\newcommand{\degree}{{\scriptscriptstyle \circ }}
\newcommand{\dfn}{\stackrel{\triangle}{=}}
\def\complex{\mathop{\raise .45ex\hbox{${\bf\scriptstyle{|}}$}
     \kern -0.40em {\rm \textstyle{C}}}\nolimits}
\def\hilbert{\mathop{\raise .21ex\hbox{$\bigcirc$}}\kern -1.005em {\rm\textstyle{H}}} 
\newcommand{\RAISE}{{\:\raisebox{.6ex}{$\scriptstyle{>}$}\raisebox{-.3ex}
           {$\scriptstyle{\!\!\!\!\!<}\:$}}} 

\newcommand{\hh}{{\:\raisebox{1.8ex}{$\scriptstyle{\degree}$}\raisebox{.0ex}
           {$\textstyle{\!\!\!\! H}$}}}

\newcommand{\OO}{\won}
\newcommand{\calA}{{\mathcal A}}
\newcommand{\calB}{{\mathcal B}}
\newcommand{\calC}{{\cal C}}
\newcommand{\calD}{{\cal D}}
\newcommand{\calE}{{\cal E}}
\newcommand{\calF}{{\mathcal F}}
\newcommand{\calG}{{\cal G}}
\newcommand{\calH}{{\cal H}}
\newcommand{\calK}{{\cal K}}
\newcommand{\calL}{{\mathcal L}}
\newcommand{\calM}{{\mathcal M}}
\newcommand{\calO}{{\cal O}}
\newcommand{\calP}{{\cal P}}
\newcommand{\calT}{{\mathcal T}} 
\newcommand{\calU}{{\mathcal U}}
\newcommand{\calX}{{\cal X}}
\newcommand{\calY}{{\mathcal Y}}
\newcommand{\calZ}{{\mathcal Z}}
\newcommand{\calXX}{{\cal X\mbox{\raisebox{.3ex}{$\!\!\!\!\!-$}}}}
\newcommand{\calXXX}{{\cal X\!\!\!\!\!-}}
\newcommand{\gi}{{\raisebox{.0ex}{$\scriptscriptstyle{\cal X}$}
\raisebox{.1ex} {$\scriptstyle{\!\!\!\!-}\:$}}}
\newcommand{\intsim}{\int_0^1\!\!\!\!\!\!\!\!\!\sim}
\newcommand{\intsimt}{\int_0^t\!\!\!\!\!\!\!\!\!\sim}
\newcommand{\pp}{{\partial}}
\newcommand{\al}{{\alpha}}
\newcommand{\sB}{{\cal B}}
\newcommand{\sL}{{\cal L}}
\newcommand{\sF}{{\cal F}}
\newcommand{\sE}{{\cal E}}
\newcommand{\sX}{{\cal X}}
\newcommand{\R}{{\rm I\!R}}
\renewcommand{\L}{{\rm I\!L}}
\newcommand{\vp}{\varphi}
\newcommand{\N}{{\rm I\!N}}
\def\ooo{\lip}
\def\ccc{\rip}
\newcommand{\ot}{\hat\otimes}
\newcommand{\rP}{{\Bbb P}}
\newcommand{\bfcdot}{{\mbox{\boldmath$\cdot$}}}

\renewcommand{\varrho}{{\ell}}
\newcommand{\dett}{{\textstyle{\det_2}}}
\newcommand{\sign}{{\mbox{\rm sign}}}
\newcommand{\TE}{{\rm TE}}
\newcommand{\TA}{{\rm TA}}
\newcommand{\E}{{\rm E\,}}
\newcommand{\won}{{\mbox{\bf 1}}}
\newcommand{\Lebn}{{\rm Leb}_n}
\newcommand{\Prob}{{\rm Prob\,}}
\newcommand{\sinc}{{\rm sinc\,}}
\newcommand{\ctg}{{\rm ctg\,}}
\newcommand{\loc}{{\rm loc}}
\newcommand{\trace}{{\,\,\rm trace\,\,}}
\newcommand{\Dom}{{\rm Dom}}
\newcommand{\ifff}{\mbox{\ if and only if\ }}
\newcommand{\nproof}{\noindent {\bf Proof:\ }}
\newcommand{\remark}{\noindent {\bf Remark:\ }}
\newcommand{\remarks}{\noindent {\bf Remarks:\ }}
\newcommand{\note}{\noindent {\bf Note:\ }}

\newcommand{\boldx}{{\bf x}}
\newcommand{\boldX}{{\bf X}}
\newcommand{\boldy}{{\bf y}}
\newcommand{\boldR}{{\bf R}}
\newcommand{\uux}{\uu{x}}
\newcommand{\uuY}{\uu{Y}}

\newcommand{\limn}{\lim_{n \rightarrow \infty}}
\newcommand{\limN}{\lim_{N \rightarrow \infty}}
\newcommand{\limr}{\lim_{r \rightarrow \infty}}
\newcommand{\limd}{\lim_{\delta \rightarrow \infty}}
\newcommand{\limM}{\lim_{M \rightarrow \infty}}
\newcommand{\limsupn}{\limsup_{n \rightarrow \infty}}

\newcommand{\ra}{ \rightarrow }

\newcommand{\ARROW}[1]
  {\begin{array}[t]{c}  \longrightarrow \\[-0.2cm] \textstyle{#1} \end{array} }

\newcommand{\AR}
 {\begin{array}[t]{c}
  \longrightarrow \\[-0.3cm]
  \scriptstyle {n\rightarrow \infty}
  \end{array}}

\newcommand{\pile}[2]
  {\left( \begin{array}{c}  {#1}\\[-0.2cm] {#2} \end{array} \right) }

\newcommand{\floor}[1]{\left\lfloor #1 \right\rfloor}

\newcommand{\mmbox}[1]{\mbox{\scriptsize{#1}}}

\newcommand{\ffrac}[2]
  {\left( \frac{#1}{#2} \right)}

\newcommand{\one}{\frac{1}{n}\:}
\newcommand{\half}{\frac{1}{2}\:}

\def\le{\leq}
\def\ge{\geq}
\def\lt{<}
\def\gt{>}

\def\squarebox#1{\hbox to #1{\hfill\vbox to #1{\vfill}}}
\newcommand{\nqed}{\hspace*{\fill}
          \vbox{\hrule\hbox{\vrule\squarebox{.667em}\vrule}\hrule}\bigskip}

\title{Variational calculation of Laplace transforms via  entropy on
   Wiener space and some Applications}

\author{Ali  S\"uleyman  \"Ust\"unel}

\begin{abstract}
\noindent
Let $(W,H,\mu)$ be the classical Wiener space where $H$ is the
Cameron-Martin space which consists of the primitives of the elements
of $L^2([0,1],\,dt)\otimes \R^d$, we denote by $L^2_a(\mu,H)$ the
equivalence classes w.r.t. $dt\times d\mu$ whose Lebesgue
densities $s\to\dot{u}(s,w)$ are almost surely adapted to the canonical
Brownian filtration. If $f$ is a Wiener functional
s.t. $\frac{1}{E[e^{-f}]}e^{-f}d\mu$ is of finite relative entropy
w.r.t. $\mu$, we prove that 
\beaa
J_\star&=& \inf\left(E_\mu\left[f\circ U+\half |u|_H^2\right]: u\in
   L_a^2(\mu,H)\right)\\
&\geq&-\log E_\mu[e^{-f}]=\inf\left(\int_W fd\ga+H(\ga|\mu):\,\nu\in
  P(W)\right)
\eeaa
where $P(W)$ is the set of probability measures on $(W,\calB(W))$ and
$H(\ga|\mu)$ is the relative entropy of $\ga$ w.r.t. $\mu$.
We call $f$ a tamed functional if the inequality above can be replaced
with equality, we characterize the class of tamed functionals, which
is much larger than the set of essentially bounded Wiener functionals.
We show that for a tamed functional  the minimization problem of l.h.s. has a  solution $u_0$ if and only if $U_0=I_W+u_0$ is almost surely
invertible and 
$$
\frac{dU_0\mu}{d\mu}=\frac{e^{-f}}{E_\mu[e^{-f}]}
$$
and then $u_0$ is unique. To do this is we prove the theorem which
says  that the relative entropy of
$U_0\mu$ is equal to the energy of $u_0$ if and only if it has a
$\mu$-a.s. left inverse. We use these results to
prove the strong  existence of the solutions of stochastic
differentail equations with singular (functional) drifts and also to
prove the non-existence of strong solutions of some stochastic
differential equations.

\end{abstract}
\maketitle
\tableofcontents
\section{\bf{Introduction}}
\noindent
Let  $(W,H,\mu)$ be  the classical Wiener space, i.e.,
$W=C_0([0,1],\R^d)$, $H$ is the corresponding Cameron-Martin space
consisting of $\R^d$-valued  absolutely continuous functions
on $[0,1]$ with square integrable derivatives w.r.t the Lebesgue
measure. Denote by  $(\calF_t,\,t\in
[0,1])$  the filtration of the canonical Wiener process, completed
w.r.t. $\mu$-negligeable sets.  Let $V:W\to W$ be a
mapping of the form $V=I_W+v$, $v:W\to H$, i.e., 
$$
V_t(w)=w(t)+v(t,w)=W_t(w)+\int_0^t\dot{v}_s(w)ds\,,
$$
where $w\to\dot{v}_s(w)$ is $\calF_s$-measurable $ds$-a.s., and
$(s,w)\to\dot{v}_s(w)$ is measurable w.r.t. the product sigma algebra
$\calB([0,1])\otimes\calF$. Heuristically, the existence of the strong
solution of the following stochastic differential equation:
$$
dU_t=-\dot{v}_t\circ U \,dt+dW_t
$$
can be interpreted as the existence of an optimal element of the
following minimization problem:
$$
K_\star=\inf\left(\half E\left[|v\circ
    (I_W+\xi)+\xi|_H^2\right]:\,\xi\in L^2_a(\mu,H)\right)
$$
where $ L^2_a(\mu,H)$ is the set of functionals as $v$ described above
with square integrable $H$-norm. The difficulty in this method lies in
the fact that, due to the quadratic character of the cost function,
the classical variational approach requires very strong regularity
hypothesis about the vector field $v$, which make the things
unrealistic. Let us write this problem in a different form: assume
that the Girsanov exponential of the vector field $v$, denoted as
$\rho(-\delta v)$ is a probability density, i.e., $E[\rho(-\delta
v)]=1$, let $f=-\log \rho(-\delta v)$. If $U=I_W+u$, with $u\in
L^2_a(\mu,H)$, we have
$$
f\circ U=\int_0^1\dot{v}_s\circ U\,dU_s+\half|v\circ U|_H^2\,.
$$
If $v\circ U\in L^2_a(\mu,H)$, taking the expectation of both sides,
we get
$$
E\left[f\circ U+\half |u|_H^2\right]=\half E\left[|v\circ
  U+u|_H^2\right]\,.
$$
Hence, again heuristically, the minimization problem $K_\star$ should
be equivalent to the minimization problem 
\beaa
J_\star&=&\inf\left(E\left[f\circ (I_W+u)+\half |u|_H^2\right]:\,u\in L^2_a(\mu,H)\right)\\
&=&\inf(J(u):\,u\in L^2_a(\mu,H))\,.
\eeaa
In this latter formulation there is no more the quadratic term
provided that the function $f$ is given directly as it happens quite
often in physics, in optimization, in the
calculation of Laplace transforms, in large deviations theory, etc. Of
course, if one studies this problem, he has to verify that the optimal
solution, if there is any, corresponds to the solution of the
corresponding stochastic differential equation. This is precisely
what we do in this paper by establishing the entropic characterization of the
$\mu$-almost sure left invertibility of the adapted perturbations of identity. 
Let us explain in more general terms the premises of the problem:
assume that  $f:W\to\R$ is  a measurable function
such that the relative  entropy of the measure
$d\nu=e^{-f}(E_\mu[e^{-f}])^{-1}d\mu$ w.r.t. $\mu$ is finite. Then it
is easy to show the validity of the following expression:
$$
-\log E[e^{-f}]=\inf\left(\int_W fd\ga+H(\ga|\mu):\,\ga\in P(W)\right)
$$
and the measure $\nu$ is the unique minimiser of the right hand side
of this equality. In case $f$ is bounded, it has been shown in \cite{B-D}
that
$$
\inf\left(E\left[f\circ (I_W+\xi)+\half |\xi|_H^2\right]:\,\xi\in L^2_a(\mu,H)\right)=-\log E[e^{-f}]\,.
$$
In case this relation holds for those $f$ which one may encounter in
the problems of invertibility mentioned above,  any minimizer $u$
for $J_\star$ will very likely  have the property that $U\mu=(I_W+u)\mu=\nu$, where
$(I_W+u)\mu$ means the push forward of the measure $\mu$ by the map
$U=I_W+u$. At this point an important concept comes up, namely we shall
call $f$ a {\bf{tamed functional}} if one has the following identity:
\begin{equation}
\label{main-def}
\inf\left(\int_W fd\ga+H(\ga|\mu):\,\ga\in
  P(W)\right)=\inf\left(E\left[f\circ (I_W+\xi+\half |\xi|_H^2\right]:\,\xi\in
  L^2_a(\mu,H)\right)\,.
\end{equation}
We prove in Theorem \ref{main-1-thm} that if $f\in L^{1+\eps}(\mu)$
for some $\eps>0$, is
such that the corresponding measure is of finite relative entropy
w.r.t. $\mu$, then it is a tamed functional. To prove the equality
between the minimizing measure of the left side of (\ref{main-def})
and the $(I_W+u)\mu$, where $u$ is the minimizing vector of the right
hand side of (\ref{main-def}), if there is any such element, 
we shall need the extension of the  results of \cite{ASU-2} 
as well as some of Ph. D. Thesis of R. Lassalle (cf.\cite{RL}). Namely, the following result will be
essential in the sequel:{\bf{ Assume that $U=I_W+u$ is an API, then it is
    $\mu$-a.s. left invertible if and only if the following equality
    holds true:}}
$$
\half E[|u|_H^2]=H(U\mu|\mu)\,,
$$
\noindent
where $H(U\mu|\mu)$ is the relative entropy of $U\mu=(I_W+u)\mu$ (push-forward of
$\mu$ under $U$) w.r.t $\mu$. Using this result we prove on the one
hand  the equivalence
between the existence of a minimizer $u$ and the $\mu$-almost sure
invertibility of the corresponding API, namely, of $U=I_W+u$ (which is
neccessarily unique) and on the other hand that $u\in L_a^2(\mu,H)$ is
a minimizing element for $J_\star$ if and only if the measure
$U\mu=(I_W+u)\mu$ is the unique minimizer of $K_\star$.  It is a remarkable fact that using this method
we can solve stochastic differential equations with very singular
(functional) drifts, e.g., we can show that the following equation is
well-defined and has a unique strong solution
$$
X_t=W_t-\int_0^t\left(\frac{E[D_\tau e^{-f}|\calF_\tau]}{E[e^{-f}|\calF_\tau]}\right)\circ
X d\tau\,,
$$
where $f$ is any $1$-convex tamed functional of $W$ and where $D_\tau F$ denotes the
density of the Sobolev derivative of $F:W\to \R$ w.r.t. the Lebesgue
measure of $[0,1]$. Note that it is not even evident to justify
$D_\tau e^{-f}$ using the classical Malliavin calculus since this
derivative may  exist only in the sense of distributions. Of course the
next step is to characterize the class of functions $f$ for which the
minimization problem is well-defined. This requires a version of
calculus of variations on the space of adapted, $dt\times d\mu$-square
integrable processes combined with the Malliavin calculus  and  with
the notion  of $H$-convex
functions;  a concept which is specific to the Wiener space(s). 

As a final application of these results, we prove that, for a given
$H$-convex subset $A\subset W$ with $\mu(A)\in (0,1)$, there is no API of the form $U=I_W+u$
which is $\mu$-a.s. left invertible and that
$dU\mu=d\nu=\frac{1}{\mu(A)}1_Ad\mu$. To understand the meaning of this
result, let us write, via It\^o representation theorem,  
$$
\frac{1_A}{\mu(A)}=\rho(-\delta v)
$$
where $v\in L^2_a(\nu,H)$. Then, under $\nu$, $V=I_W+v$ is a Brownian
motion and  the SDE
$$
dU_t=-\dot{v}_t\circ U dt+dV_t
$$
has a weak solution but has no strong solution. As the reader can
realize,  this
result is a consequence of the variational calculus developed here and
hence its nature and its philosophy are quite different from the example of
B. Tsirelson, cf. \cite{Tsi}.

Let us finally add that the results of this paper have immediate
extensions  to the infinite dimensional case (i.e., the cylindrical
Brownian motion) and also to 
the abstract Wiener spaces via the theory developed in \cite{FILT}.
\section{\bf{Preliminaries and notation}}
\label{preliminaries}
Let $W$ be the classical Wiener  space with  the Wiener
measure $\mu$. The
corresponding Cameron-Martin space is denoted by $H$. Recall that the
injection $H\hookrightarrow W$ is compact and its adjoint is the
natural injection $W^\star\hookrightarrow H^\star\subset
L^2(\mu)$. A subspace $F$ of $H$ is called regular if the
corresponding orthogonal projection
has a continuous extension to $W$, denoted again  by the same letter.
It is well-known that there exists an increasing sequence of regular
subspaces $(F_n,n\geq 1)$, called total,  such that $\cup_nF_n$ is
dense in $H$ and in $W$. Let $\sigma(\pi_{F_n})${\footnote{For the notational
  simplicity, in the sequel we shall denote  it by  $\pi_{n}$.}}  be the
$\sigma$-algebra generated by $\pi_{F_n}$, then  for any  $f\in
L^p(\mu)$, the martingale  sequence
$(E[f|\sigma(\pi_{F_n})],n\geq 1)$
converges to $f$ (strongly if  $p<\infty$) in $L^p(\mu)$. Observe that the function
$f_n=E[f|\sigma(\pi_{F_n})]$ can be identified with a function on the
finite dimensional abstract Wiener space $(F_n,\mu_n,F_n)$, where
$\mu_n=\pi_n\mu$.

Since the translations of $\mu$ with the elements of $H$ induce measures
equivalent to $\mu$, the G\^ateaux  derivative in $H$ direction of the
random variables is a closable operator on $L^p(\mu)$-spaces and  this
closure will be denoted by $\nabla$ cf.,  for example \cite{ASU, ASU-1}. The corresponding Sobolev spaces
(the equivalence classes) of the  real random variables
will be denoted as $\DD_{p,k}$, where $k\in \NN$ is the order of
differentiability and $p>1$ is the order of integrability. If the
random variables are with values in some separable Hilbert space, say
$\Phi$, then we shall define similarly the corresponding Sobolev
spaces and they are denoted as $\DD_{p,k}(\Phi)$, $p>1,\,k\in
\NN$. Since $\nabla:\DD_{p,k}\to\DD_{p,k-1}(H)$ is a continuous and
linear operator its adjoint is a well-defined operator which we
represent by $\delta$.  $\delta$ coincides with the It\^o
integral of the Lebesgue density of the adapted elements of
$\DD_{p,k}(H)$ (cf.\cite{ASU,ASU-1}).

For any $t\geq 0$ and measurable $f:W\to \reals_+$, we note by
$$
P_tf(x)=\int_Wf\left(e^{-t}x+\sqrt{1-e^{-2t}}y\right)\mu(dy)\,,
$$
it is well-known that $(P_t,t\in \reals_+)$ is a hypercontractive
semigroup on $L^p(\mu),p>1$,  which is called the Ornstein-Uhlenbeck
semigroup (cf.\cite{ASU,ASU-1}). Its infinitesimal generator is denoted
by $-\calL$ and we call $\calL$ the Ornstein-Uhlenbeck operator
(sometimes called the number operator by the physicists). The
norms defined by
\begin{equation}
\label{norm}
\|\phi\|_{p,k}=\|(I+\calL)^{k/2}\phi\|_{L^p(\mu)}
\end{equation}
are equivalent to the norms defined by the iterates of the  Sobolev
derivative $\nabla$. This observation permits us to identify the duals
of the space $\DD_{p,k}(\Phi);p>1,\,k\in\NN$ by $\DD_{q,-k}(\Phi')$,
with $q^{-1}=1-p^{-1}$,
where the latter  space is defined by replacing $k$ in (\ref{norm}) by
$-k$, this gives us the distribution spaces on the Wiener space $W$
(in fact we can take as $k$ any real number). An easy calculation
shows that, formally, $\delta\circ \nabla=\calL$, and this permits us
to extend the  divergence and the derivative  operators to the
distributions as linear,  continuous operators. In fact
$\delta:\DD_{q,k}(H\otimes \Phi)\to \DD_{q,k-1}(\Phi)$ and
$\nabla:\DD_{q,k}(\Phi)\to\DD_{q,k-1}(H\otimes \Phi)$ continuously, for
any $q>1$ and $k\in \reals$, where $H\otimes \Phi$ denotes the
completed Hilbert-Schmidt tensor product (cf., for instance
\cite{ASU,ASU-1}). Finally, in the case of classical Wiener space,
we denote by $\DD_{p,k}^a(H)$ the subspace defined by 
$$
\DD_{p,k}^a(H)=\{\xi\in\DD_{p,k}(H):\,\dot{\xi}\mbox{ is adapted}\}
$$
for $p\geq 1,\,k\in\R$, for $p=2,\,k=0$, we shall write $L^2_a(\mu,H)$.

\noindent
A measurable  function
$f:W\to \reals\cup\{\infty\}$  is called $\alpha$-convex, $\alpha\in\R$, if the map
$$
h\to f(x+h)+\frac{\alpha}{2}|h|_H^2=F(x,h)
$$
is convex on the Cameron-Martin space $H$ with values in
$L^0(\mu)$. Note that this notion is compatible with the
$\mu$-equivalence classes of random variables thanks to the
Cameron-Martin theorem. It is proven in \cite{F-U1} that
this definition  is equivalent  the following condition:
  Let $(\pi_n,n\geq 1)$ be a sequence of regular, finite dimensional,
  orthogonal projections of
  $H$,  increasing to the identity map
  $I_H$. Denote also  by $\pi_n$ its  continuous extension  to $W$ and
  define $\pi_n^\bot=I_W-\pi_n$. For $x\in W$, let $x_n=\pi_nx$ and
  $x_n^\bot=\pi_n^\bot x$.   Then $f$ is $1$-convex if and only if
$$
x_n\to \frac{1}{2}|x_n|_H^2+f(x_n+x_n^\bot)
$$
is  $\pi_n^\bot\mu$-almost surely convex.

\noindent
We shall use also the following result, which makes part of the
folklore of the Wiener measure:
\begin{lemma}
\label{exp-lemma}
Denote by $\rho(\delta h),\, h\in H$, the Wick exponential
$$
\rho(\delta h)=\exp\left(\delta h-\half|h|_H^2\right)
$$
where $\delta h=\int_0^1\dot{h}(s)dW_S$ (i.e., Wiener integral). The
map
$$
h\to \rho(\delta h)
$$
is weakly continuous on $H$ with values in $L^p(\mu)$, for any $p>1$.
\end{lemma}
\nproof
Assume that $F$ is of finite Wiener chaos, then, for any $h\in H$, it
follows from Cameron-Martin theorem that 
\beaa
E[F\rho(\delta h)]&=&E[F(\cdot+h)]\\
&=&\sum_{n=0}^\infty\frac{1}{n!}\left(E[\nabla^nF, h^{\otimes
    n}\right)_{H^{\otimes n}}
\eeaa
(cf. \cite{ASU,ASU-1}).  If $h_k\to h$ weakly in $H$, then
$h_k^{\otimes m}\to h^{\otimes m}$
weakly in $H^{\otimes m}$ for any $m\geq 1$, hence $E[F\rho(\delta
h_k)]\to E[F\rho(\delta h)]$ as $k\to\infty$ if $F$ is chosen as
above. If $F\in L^p(\mu)$, then there exists a sequence $(F_n,n\geq
1)$ converging to $F$ in $L^p(\mu)$ with each $F_n$ being of finite
Wiener chaos. Let $c_n(h)=E[F_n\rho(\delta h)],\,c(h)=E[F\rho(\delta
h)]$, then 
$$
|c(h)-c(k)|\leq|c(h)-c_n(h)|+|c_n(h)-c_n(k)|+|c_n(k)-c(k)|
$$
hence $\lim c(h)=c(k)$ as $h\to k$ weakly in $H$.
\nqed

\begin{definition}
A map $u:W\to H$ is called an $H-C^1$-map if the map $h\to u(w+h)$ is
Fr\'echet differentiable on $H$ for almost all $w$.
\end{definition}
\begin{remarkk}
\label{rem-1}
This is a very strong property, in particular it implies that the set
of the elements  $w$ of $W$ should be $H$-invariant. Let us note that,
if $u$ is in some space $L^p(\mu,H)$, then $P_\tau u$ is an
$H-C^1$-map for any $\tau>0$, where $P_\tau$ denotes the
Ornstein-Uhlenbeck semigroup on $W$ (in fact it is even $H$-analytic,
cf.\cite{BOOK}, Chapter 2).
\end{remarkk}

\noindent
A more relaxed notion is given as
\begin{definition}
\label{loc}
The map $u$ is called an $(H-C^1)_{\rm loc}$ map if there exists an
almost surely strictly positive map $q$ such that $h\to u(w+h)$ is
continuously differentiable on the set $\{h\in H:\,|h|_H<q(w)\}$.
\end{definition}
We have the following result about the change of variables formula for
$(H-C^1)_{\rm loc}$-maps proved in Theorem 4.4.1 of Chapter 4 of
\cite{BOOK} 
\begin{theorem}
\label{jacobi-thm}
Assume that $u\in L^2_a(\mu,H)$ is an $(H-C^1)_{\rm loc}$-mapping,
define $U:W\to W$ as $U=I_W+u$,  let
$Q$ be the set $\{w\in W: q(w)>0\}$, where $q$ is the mapping given in
Definition \ref{loc}, then for any $f,g\in C_B(W)$, the following
identity holds true:
$$
E[f\circ U\,g\,\rho(-\delta u)]=E\left[f\sum_{y\in U^{-1}\{w\}\cap
    Q}g(y)\right]\,.
$$
\end{theorem}

\section{\bf{Characterization of the invertible shifts}}
\noindent
We  begin with the definition of the notion of  almost sure
invertibility with respect to a measure. This notion is extremely
important since it makes the  things work

\begin{definition}
\begin{itemize}
\item A measurable map $T:W\to W$ is called ($\mu$-) almost surely left invertible
if there exists a measurable  map $S:W\to W$ such that 
and $S\circ T =I_W$ $\mu$-a.s.
\item  Moreover, in this case it is trivial to
see that $T\circ S=I_W$ $T\mu$-a.s., where $T\mu$ denotes the image of
the measure $\mu$ under the map $T$. 
\item If $T\mu$ is equivalent to $\mu$,
then we say in short that $T$ is $\mu$-a.s. invertible. 
\item Otherwise, we
may say that $T$ is $(\mu,T\mu)$-invertible in case  precision is
required or just $\mu$-a.s. left invertible and $S$ is called the
$\mu$-left inverse of $T$.
\end{itemize}
\end{definition}

\begin{theorem}
\label{thm-0}
For any $u\in L_a^2(\mu,H)$, we have the following inequality
$$
H(U\mu|\mu)\leq \frac{1}{2}E\int_0^1|\dot{u}_s|^2ds\,,
$$
where $H(U\mu|\mu)$ is the relative entropy of the measure $U\mu$
w.r.t. $\mu$. 
\end{theorem}
\nproof
Let $L$ be the Radon-Nikodym density of $U\mu$ w.r.t. $\mu$. For any
$0\leq g\in C_b(W)$, using the Girsanov theorem, we have 
$$
E[g\circ U]=E[g\,L]\geq E[g\circ U\,L\circ U\,\rho(-\delta u)]\,,
$$
hence 
$$
L\circ U\,E[\rho(-\delta u)|U]\leq 1
$$
$\mu$-a.s. Consequently, using the Jensen inequality
\beaa
H(U\mu|\mu)&=&E[L\log L]=E[\log L\circ U]\\
&\leq&-E[\log E[\rho(-\delta u)|U]]\\
&\leq&-E[\log \rho(-\delta u)]\\
&=&\half E\int_0^1|\dot{u}_s|^2ds\,.
\eeaa
\nqed

\begin{theorem}
\label{main-thm}
Assume that $U=I_W+u$ is an API, i.e., $u\in L^2_a(\mu,H)$ such that
$s\to\dot{u}(s,w)$ is $\calF_s$-measurable for almost all $s$. 
Then $U$ is almost surely left invertible with a left inverse $V$ if and only if 
$$
H(U\mu|\mu)=\frac{1}{2}E[|u|_H^2]=\frac{1}{2}E\int_0^1|\dot{u}_s|^2ds\,,
$$
i.e., if and only if the entropy of $U\mu$ is equal to the energy of
the drift $u$.
\end{theorem}
\nproof
Due to Theorem \ref{thm-0}, the relative entropy is finite as soon as 
$u\in L^2_a(\mu,H)$. Let us suppose now that the equality holds and let us denote by $L$ the Radon-Nikodym derivative of
$U\mu$ w.r.t. $\mu$. Using the It\^o representation theorem, we can
write 
$$
L=\exp\left(-\int_0^1 \dot{v}_sdW_s-\half\int_0^1|\dot{v}_s|^2ds\right)
$$
$U\mu$-almost surely. Let $V=I_W+v$, as described in \cite{Foll}, from the It\^o formula and Paul
L\'evy's theorem, it is immediate that $V$ is an $U\mu$-Wiener
process, hence 
\begin{equation}
\label{ent}
E[L\log L]=\half E[L\,|v|_H^2]\,.
\end{equation}
Now, for any $f\in C_b(W)$, we have from the Girsanov
theorem
$$
E[f\circ U]=E[f\,L]\geq E[f\circ U\,L\circ U\,\rho(-\delta u)]
$$
consequently
$$
L\circ U\,E[\rho(-\delta u)|U]\leq 1
$$
$\mu$-a.s. Let us denote $E[\rho(-\delta u)|U]$ by $\hat{\rho}$. We
have then $\log L\circ U+\log\hat{\rho}\leq 0$ $\mu$-a.s. Taking the
expectation w.r.t. $\mu$ and the Jensen inequality give
\beaa
H(U\mu|\mu)&=&E[L\log L]\leq -E[\log \hat{\rho}]\\
&\leq&-E[\log\rho(-\delta u)]=\half E[|u|_H^2]\,.
\eeaa
Since $\log$ is a strictly concave function, the equality $E[\log
\hat{\rho}]=E[\log\rho(-\delta u)]$ implies that $\rho(-\delta u)=\hat{\rho}$
$\mu$-a.s. Hence we obtain
\begin{equation}
\label{eqq}
E[L\log L+\log\rho(-\delta u)]=E[\log(L\circ U\,\rho(-\delta u))]=0\,,
\end{equation}
since $L\circ U\rho(-\delta u)\leq 1$ $\mu$-a.s., the equation
(\ref{eqq}) implies
\begin{equation}
\label{basic}
L\circ U\rho(-\delta u)=1
\end{equation}
$\mu$-a.s. Combining  the exponential representation of $L$ with the
relation (\ref{basic}) implies
\bea
\label{basic-1}
0&=&\left(\int_0^1\dot{v}_sdW_s\right)\circ U+\half|v\circ U|_H^2+\delta
u+\half|u|_H^2\nonumber\\
&=&\delta(v\circ U)+\delta u+(v\circ U,u)_H+\half(|u|_H^2+|v\circ
U|_H^2)\nonumber\\
&=&\delta(v\circ U+u)+\half|v\circ U+u|_H^2
\eea
$\mu$-a.s. From the relation (\ref{ent}) it follows that $v\circ U\in
L^2_a(\mu,H)$, hence taking the expectations of both sides of
(\ref{basic-1}) w.r.t. $\mu$ is licit  and this implies $v\circ U+u=0$ $\mu$-a.s., which means that
$V=I_W+v$ is the $\mu$-left inverse of $U$. 

To show the neccessity, let us denote by $(L_t,t\in [0,1])$ the
martingale
$$
L_t=E[L|\calF_t]=E\left[\frac{dU\mu}{d\mu}|\calF_t\right]
$$
and let 
$$
T_n=\inf\left(t:\,L_t<\frac{1}{n}\right)\,.
$$
Since $U\circ V=I_W$ $(U\mu)$a.s.,  $V$ can be written as $V=I_W+v$
($U\mu$)-a.s. and that $v\in L_a^0(U\mu,H)$, i.e., 
$v(t,w)=\int_0^t \dot{v}_s(w)ds$, $\dot{v}$
is adapted to the filtration $(\calF_t)$ completed w.r. to $U\mu$ and
$\int_0^1|\dot{v}_s|^2ds<\infty$ $(U\mu)$-a.s. Since $\{t\leq
T_n\}\subset \{L>0\}$ and since on this latter set $\mu$ and $U\mu$
are equivalent, we have 
$$
\int_0^{T_n}|\dot{v}_s|^2ds<\infty
$$
{\bf{$\mu$-almost surely}}. Consequently the  inequality 
$$
E_\mu[\rho(-\delta v^n)]\leq 1
$$
holds true for any $n\geq 1$, where $v^n(t,w)=\int_0^t
1_{[0,T_n]}(s,w)\dot{v}_s(w)ds$. By positivity we also have 
$$
E_\mu[\rho(-\delta v^n)1_{\{L>0\}}]\leq 1\,.
$$
Since $\lim_nT_n=\infty $$(U\mu)$-a.s., we also have $\lim_nT_n=\infty
$ $\mu$-a.s. on the set $\{L>0\}$ and the Fatou lemma implies
\begin{equation}
\label{ineq-1}
E_\mu[\rho(-\delta v)1_{\{L>0\}}] =E_\mu[\lim_n\rho(-\delta
v^n)1_{\{L>0\}}]\leq \lim\inf_n E_\mu[\rho(-\delta v^n)1_{\{L>0\}}]\leq 1\,.
\end{equation}
for any $n\geq 1$.  From the identity $U\circ V=I_W$
$(U\mu)$-a.s., we have $v+u\circ V=0$ $(U\mu)$-a.s., hence $v\circ
U+u=0$ $\mu$-a.s.  An algebraic  calculation gives immediately 
\begin{equation}
\label{ineq-2}
\rho(-\delta v)\circ U\,\rho(-\delta u)=1
\end{equation}
$\mu$-a.s. Now applying the Girsanov theorem to API $U$ and using the
relation (\ref{ineq-2}), we obtain 
\beaa
E[g\circ U]&=&E[g\,L]=E\left[g\circ U(\rho(-\delta v)1_{\{L>0\}})\circ
U\rho(-\delta u)\right]\\
&\leq&E\left[g\,\rho(-\delta v)1_{\{L>0\}}\right]\,,
\eeaa
for any positive $g\in C_b(W)$ (note that on the set $\{L>0\}$
$\rho(-\delta v)$ is perfectly well-defined w.r. to $\mu$). Therefore 
$$
L\leq \rho(-\delta v)1_{\{L>0\}}
$$
$\mu$-a.s. Now, this last inequality, combined with the inequality
(\ref{ineq-1}) entails that 
$$
L=\rho(-\delta v)1_{\{L>0\}}
$$
$\mu$-a.s., hence 
$$
L\circ U\,\rho(-\delta u)=1
$$
$\mu$-a.s. To complete the proof it suffices to remark then that 
\beaa
H(U\mu|\mu)&=&E[L\log L]=E[\log L\circ U]\\
&=&E[-\log \rho(-\delta u)]\\
&=&
\half E[|u|_H^2]\,.
\eeaa
\nqed

\noindent
The following theorem, although it has strong hypothesis, is at the
heart of the further developements:
\begin{theorem}
\label{H-C-thm}
Assume that $u\in L^2_a(\mu,H)$ is an $(H-C^1)_{\rm loc}$ such that 
$E[\rho(-\delta u)]=1$, then the mapping $U=I_W+u$ is
$\mu$-a.s. invertible.
\end{theorem}
\nproof
We have, from Theorem \ref{jacobi-thm}, taking $g=1$ and defining the
multiplicity of $U$ on the set $Q$ (cf. the notations of Theorem
\ref{jacobi-thm}) as
$$
N(w,Q)=\sum_{y\in U^{-1}\{w\}\cap Q} 1(y)\,,
$$
the relation
$$
E[f\circ U\,\rho(-\delta u)]=E[f\,N(\cdot,Q)]\,.
$$
On the other hand the Girsanov theorem implies that 
$$
 E[f\circ U\,\rho(-\delta u)]=E[f]\,,
$$
for any $f\in C_b(W)$, hence $N(w,Q)=1$ $\mu$-a.s., this implies that
the map $U$ is almost surely injective. The hypothesis $E[\rho(-\delta
u)]=1$ implies also that $U(W)=W$ $\mu$-a.s., hence $U$ is almost
surely surjective. 
\nqed
\begin{remarkk}
Again using the Girsanov theorem, it is immediate to show that 
the inverse of $U$ is of the form $V=I_W+v$, with $v\in L^0_a(\mu,H)$.
\end{remarkk}

\section{\bf{Some variational problems related to entropy and large
    deviations}}
\noindent
The following is an extension of a well-known result in large
deviations theory:
\begin{theorem}
v\label{Fenchel}
Let $(A,\calA)$ be a measurable space and let $f:A\to \R$ be a
measurable function, denote by $P(A)$ the set of probability measures
on $(A,\calA)$. Suppose that for some $\gamma\in P(A)$, $f$ satisfies
$$
\int_A(|f|(1+e^{f}))d\gamma<\infty\,.
$$
Then the following identity holds:
$$
\log\gamma(e^{f})=\sup\left(\int fd\nu-H(\nu|\gamma):\,\nu\in P(A)\right)
$$
and the unique supremum is attained at the measure
$$
d\nu=\frac{e^{f}}{\gamma[e^{f}]}d\gamma\,,
$$
where $H(\nu|\gamma)$
denotes the relative entropy of $\nu$ w.r.t. $\gamma$.
\end{theorem}
\nproof
When $f$ is bounded, this theorem is well-known (cf.\cite{B-D} and the
references there). First, suppose that $f$ is lower-bounded, let
$f_n=f\wedge n$, it follows from the bounded case that 
$$
\log\int e^fd\gamma\geq \sup\left(\int fd\nu-H(\nu|\ga)\right)\,.
$$
Since $\nu(f_n)\leq \nu(f)$ for any $\nu$ with $H(\nu|\ga)<\infty$, we
also have 
$$
\log \ga(e^{f_n})=\sup_\nu\left(\int f_nd\nu-H(\nu|\ga)\right)\leq \sup_\nu\left(\int
f d\nu-H(\nu|\ga)\right)
$$
and passing to the limit we get 
$$
\log \ga(e^{f})\leq \sup_\nu\left(\int f d\nu-H(\nu|\ga)\right)
$$
and this proves the claim when $f$ is lower bounded.
For the general case, define $g_\eps=\log(e^f+\eps)$, then $g_\eps$ is
lower bounded, hence
$$
\log\ga(e^{g_\eps})\geq \int g_\eps d\nu-H(\nu|\ga)
$$
for any $\eps$ and for any $\nu$ with finite $H(\nu|\ga)$. Passing to
the limit as $\eps\to 0$ and taking the supremum w.r.t. $\nu$,  we get 
$$
\log\ga(e^{f})\geq \sup_\nu\left(\int f d\nu-H(\nu|\ga)\right)\,.
$$
To see the equality, it suffices to remark that for the measure 
$$
d\nu_0=\frac{e^f}{\ga(e^f)}d\ga
$$
the supremum is attained.
\nqed

\begin{remarkk}
In the sequel, we shall use a variation of this theorem where $f$ will
be replaced by $-f$ and the corresponding equality is 
$$
-\log\int e^{-f}d\ga=\inf\left(\int fd\nu+H(\nu|\ga): \,\nu\in P(A)\right)\,.
$$
\end{remarkk}

\begin{theorem}
\label{1st-ineq}
Assume that 
$$
\int_W(|f|+1)e^{-f}d\mu<\infty\,.
$$
Then the following inequality holds true:
$$
-\log E[e^{-f}]\leq \inf\left(E\left[f\circ (I_W+u)+\half|u|_H^2\right]:\,u\in
  L^2(\mu,H)\right)\,.
$$
\end{theorem}
\nproof
Combining Theorem \ref{Fenchel} with Remark 3, we see already that $f$
is quasi-integrable for any measure $\nu$ which is of finite relative
entropy w.r.t. the Wiener measure $\mu$. In particular, if
$\nu=(I_W+u)(\mu)$, with $u\in L^2_a(\mu,H)$, then, from Theorem
\ref{thm-0}, we have $H((I_W+u)\mu|\mu)\leq \half\|u\|^2_{L^2(\mu,H)}$,
  hence the inequality follows.
\nqed


\begin{theorem}
\label{main-1-thm}
Assume that $f\in L^p(\mu),\,e^{-f}\in L^q(\mu)$  with
$p^{-1}+q^{-1}=1$, hence
$$
\int (|f|+1)e^{-f}d\mu<\infty\,.
$$
Then the following equalities hold true
\beaa
J_\star&=&-\log\mu(e^{-f})\\
&=&\inf\left(\int_W f d\nu+H(\nu|\mu):\,\nu\in P(W)\right)\\
&=&\inf\left[\int
  [f\circ(I_W+u)+\half|u|_H^2]d\mu:\,u\in L^2_a(\mu,H)\right]
\eeaa
\end{theorem}
\nproof
We just need to prove the last equality; we shall proceed the proof by
showing that each side of this last equality is less than the other
one. First, it is immediate from the definition of infimum and from
Theorem \ref{thm-0} that 
$$
-\log\mu(e^{-f})\leq \inf\left[\int
  [f\circ(I_W+u)+\half|u|_H^2]d\mu:\,u\in L^2_a(\mu,H)\right]\,.
$$
To show the reverse inequality is more delicate: let $(e_n,n\geq 1)$
be a complete, orthonormal basis of the Cameron-Martin space $H$,
denote by $V_n,\,n\geq 1$, the sigma-algebra generated by the Gaussian
random variables $\delta e_1,\ldots,\delta e_n$. Define now $f_n$ as 
$$
f_n=E[P_{1/n}f|V_n]
$$ 
where $P_{1/n}$ is the Ornstein-Uhlenbeck semi-group on $W$. Denote by
$l_n$ the density $e^{-f_n}/E[e^{-f}]$ and define
$$
\dot{v}^n_t=\frac{E[D_te^{-f_n}|\calF_t]}{E[e^{-f_n}|\calF_t]}
$$
where $D_t e^{-f_n}$ denotes the Lebesgue density of the
$H$-derivative $\nabla e^{-f_n}$ which is perfectly well-defined
thanks to the hypothesis. Let $v^n$ be the primitive of $\dot{v}^n$,
i.e., 
$$
v^n(t,w)=\int_0^t\dot{v}_s^n(w)ds\,.
$$
Let us indicate that the mapping
$$
w\to\int_0^\cdot E[D_te^{-f_n}|\calF_t](w) dt
$$
is an $H-C^1$-map due to the regularization with the
Ornstein-Uhlenbeck semi-group (cf. Remark \ref{rem-1}). Let $H_\sigma$
denote the Cameron-Martin space equipped with its weak topology, then,
for any $h\in H$,
\beaa
E[f_n|\calF_t](w+h)&=&\int_W
E[E[f|V_n]|\calF_t](e^{-1/n}(w+h)+\sqrt{1-e^{-2/n}}y)\mu(dy)\\
&=&\int_w\rho(\al_n\delta
h(y))E[E[f|V_n]|\calF_t](e^{-1/n}w+\sqrt{1-e^{-2/n}}y)\mu(dy)\,,
\eeaa
where $\al_n=e^{-1/n}/\sqrt{1-e^{-2/n}}$. It follows from Lemma
\ref{exp-lemma} and from the hypothesis about $f$ that the mapping 
$$
(t,h)\to E[E[P_{1/n}f|V_n]|\calF_t](w+h)
$$
is $\mu$-a.s. continuous on the space $[0,1]\times H_\sigma$.
Consequently
$$
\sup_{t\in [0,1],\,h\in B} E[E[P_{1/n}f|V_n]|\calF_t](w+h)>0
$$
$\mu$-a.s. for any bounded, weakly closed set $B\in H$.  and the set of such $w$'s are again $H$-invariant

$$
\inf_{t\in [0,1],\,h\in B} E[E[P_{1/n}e^{-f}|V_n]|\calF_t](w+h)>0
$$
$\mu$-a.s. and the set of such $w$'s are again $H$-invariant. This
observation, combined with the $H-C^1$-property of $w\to\int_0^\cdot
E[D_te^{-f_n}|\calF_t](w) dt$ implies that $v^n$ is an $H-C^1$-map and
it follows from Theorem \ref{H-C-thm} that the mapping $w\to
w+v^n(w)=V_n(w)$ is $\mu$-a.s. invertible. Let $U_n=I_W+u^n$ be its inverse,
then clearly 
$$
\frac{dU_n\mu}{d\mu}=\rho(-\delta v^n)=l_n=\frac{e^{-f_n}}{E[e^{-f_n}]}\,.
$$
It is now trivial to see from Jensen's inequality  that $u^n\in
L^2_a(\mu,H)$. Moreover
\beaa
-\lim_n\log E[e^{-f_n}]&=&\lim_n\left(\frac{1}{E[e^{-f_n}]}\int_W f_n e^{-f_n}d\mu+\half E[|u^n|_H^2]\right)\\
&=&\lim_n\left(\frac{1}{E[e^{-f_n}]}\int_W f e^{-f_n}d\mu+\half E[|u^n|_H^2]\right)\\
&=&\lim_n\left(\int_W f\circ U_n d\mu+\half E[|u^n|_H^2]\right)\\
&\geq&\inf\left(\int_W (f\circ U+|u|_H^2)d\mu:\,u\in
  L^2_a(\mu,H)\right)
\eeaa
and this completes the proof.
\nqed

\begin{definition}
We call a measurable  map $f: W\to \R\cup\{\infty\}$ with the property
$E[(1+|f|)e^{-f}]<\infty$,  a {\bf{tamed functional}} if the conclusion
of Theorem \ref{main-1-thm} is valid for $f$, mainly if
$$
-\log E[e^{-f}]=\inf\left[\int
  [f\circ(I_W+u)+\half|u|_H^2]d\mu:\,u\in L^2_a(\mu,H)\right]
$$
\end{definition}

\noindent
An immediate consequence of the logarithmic Sobolev inequality (cf. \cite{fed}) for the
Wiener measure gives
\begin{proposition}
Assume that $f\in L^{1+\eps}(\mu)$ such that $E[e^{-f}]<\infty$. Let
$f^-=\max(-f,0)$, if $E[f^-e^{f^-}]<\infty$, then $f$ is a tamed
functional, in particular the latter condition is satisfied if
$$
E[e^{-f}|\nabla f^-|^2_H]<\infty\,.
$$
\end{proposition}

\section{\bf{Caracterization of the minimizers}}

\noindent
We come  to the minimization problem for:
$$
J_\star=-\log\mu(e^{-f})=\inf\left[\int
  (f\circ(I_W+u)+\half|u|_H^2)d\mu:\,u\in L^2_a(\mu,H)\right]
$$
The following result gives a complete characterization of the
attainability of $J_\star$ in the situation of finite entropy:

\begin{theorem}
\label{attn-thm}
Assume that $f$ is a tamed functional, 
then the infimum $J_\star$ is attained at some $u\in L^2_a(\mu,H)$ if
and only if the API defined as $U=I_W+u$ has a left inverse $V=I_W+v$
with $v\in L^2_a(U\mu,H)$ and 
$$
\frac{dU\mu}{d\mu}=\frac{e^{-f}}{E[e^{-f}]}=L=\rho(-\delta v)\,.
$$
Moreover 
$U$ is the unique strong solution of the following SDE
$$
dU_t=-\dot{v}_t\circ U dt+dW_t
$$
and if $E[e^{-(1+\eps)f}]<\infty$ for some $\eps>0$, then  $\dot{v}$ can be expressed as
$$
\dot{v}_\tau=\frac{E[D_\tau
    L|\calF_\tau]}{E[L|\calF_\tau]}
$$
$d\tau\times dU\mu$-almost surely.
\end{theorem}
\nproof
Sufficiency: since $U$ is a.s. left invertible, we have from Theorem
\ref{main-thm}
$$
H(U\mu|\mu)=H(L\cdot\mu|\mu)=\half E[|u|_H^2]\,,
$$
hence it is a trivial calculation of the entropy to see that 
$$
E[f\circ U+\half|u|_H^2]=-\log E[e^{-f}]\,,
$$
hence $J_\star=J(u)$.

To prove the neccessity, suppose that there exists some $u\in
L^2_a(\mu,H)$ with $J_\star=J(u)$. Assume that $U=I_W+u$ is not
a.s. left invertible, then from Theorem \ref{main-thm}, we have 
$$
H(U\mu|\mu)<\half E[|u|_H^2]\,.
$$
Hence
$$
J_\star=E[f\circ U+\half |u|_H^2]>E[f\circ U]+H(U\mu|\mu)\,,
$$
but $f$ is a tamed functional, hence $J_\star=K_\star$ and the last
inequality  is a contradiction to the fact that $K_\star$ is the infimum of
such expressions. Therefore $H(U\mu|\mu)$ should be
equal to the energy of $u$, i.e., to $\half E[|u|_H^2]$, which is
equivalent to the left a.s. invertibility. Since the minimizing
measure of $K_\star$ is unique, we should have evidently 
$$
\frac{dU\mu}{d\mu}=\frac{e^{-f}}{E[e^{-f}]}=L\,.
$$
The expression for $L$ is obviuous from the stochastic integral
representation of Wiener functionals which are not neccessarily
Sobolev differentiable, cf.\cite{repres}.
\nqed
\begin{remarkk}
 Notice that if $f<\infty$ $\mu$-a.s. then $U$ is
$\mu$-a.s. invertible. 
\end{remarkk}

\begin{theorem}
\label{uniqueness-thm}
Assume that $f$ is a tamed functional, if $J_\star$ is attained  at
some  $u\in L_a^2(\mu,H)$, then $u$ is unique.
\end{theorem}
\nproof
Suppose that there are two such elements of $L_a^2(\mu,H)$, say
$u_1,\,u_2$ such that $J(u_1)=J(u_2)=J_\star$. Since, from Theorem
\ref{attn-thm}, 
$$
\frac{e^{-f}}{E[e^{-f}]}=\frac{dU_1\mu}{d\mu}=\frac{dU_2\mu}{d\mu}\,,
$$
where $U_i=I_W+u_i,\,i=1,2$. Moreover, if we denote $L$ as
$\rho(-\delta v)$ $\nu$-a.s., where $\nu=U_i\mu$, we see that $V\circ
U_1=V\circ U_2$ $\mu$-a.s., where $V=I_W+v$. Consequently $U_1\circ
V=U_2\circ V$ $\nu$-a.s., and it follows then that $U_1\circ V\circ
U_1=U_2\circ V\circ U_1$ $\mu$-a.s., consequently $U_1=U_2$ $\mu$-a.s.
\nqed

\begin{theorem}
\label{var-attn-thm}
Assume that $\nu$ be a probability on $(W,\calF)$ absolutely
continuous w.r.t. $\mu$, denote by $K$ the corresponding Radon-Nikodym
derivative. Let $f:W\to \R$ be a Borel function such that $\nu(|f|\exp
-f)<\infty$. Assume that $f-\log K$ is a tamed functional.
Then we have 
\beaa
-\log\nu(e^{-f})&=&\inf\left(\int_W fd\beta+H(\beta|\nu);\,\beta\in
  P(W)\right)\\
&=&\inf\left(E_\mu\left[(f-\log K)\circ
    (I_W+u)+\half|u|_H^2\right]:\,u\in L^2_a(\mu,H)\right)
\eeaa
and the second infimum is attained if and only if $U=I_W+u$ is
$\mu$-a.s. left invertible.
\end{theorem}
\nproof
The first equality follows from Theorem \ref{Fenchel},  the second
follows from the hypothesis by noting  that
$\nu(e^{-f})=E_\mu[e^{-f+\log K}]$, hence the proof follows from
Theorem \ref{attn-thm}.
\nqed


\section{\bf{Existence for $H$- convex functionals}}

\begin{theorem}
\label{conv-1}
Assume that $f\in L^0(\mu)$ is $1$-convex and that $f^-=\max(-f,0)$ is
exponentially integrable, i.e., 
$E[\exp cf^-]<\infty$ for some $c>1$. Then $J_\star$ is attained at some $u\in L^2_a(\mu,H)$ provided that
$E[f\circ (I_W+\xi)]<\infty$ for at least one $\xi\in
L_a^2(\mu,H)$. Moreover, if $f\in L^{1+\eps}(\mu)$ for some $\eps>0$,
then $f$ is a tamed functional, consequently the conclusions of
Theorem \ref{attn-thm} hold true for $f$, in particular
$$
\frac{dU\mu}{d\mu}=\frac{e^{-f}}{E[e^{-f}]}=\rho(-\delta v)
$$
where 
$$
\dot{v}_t=\frac{E[D_te^{-f}|\calF_t]}{E[e^{-f}|\calF_t]}
$$
and $V=I_W+v$ is the $\mu$-left inverse to  $U=I_W+u$ which is the
unique strong solution of the following stochastic differential equation:
$$
dU_t=-\dot{v}_t\circ U dt+dW_t\,.
$$
Finally $U$ is also the solution of the following Monge-Amp\'ere
equation:
$$
\frac{dU\mu}{d\mu}=\exp\left(-\int_0^1E_{U\mu}[D_tf|\calF_t]dW_t-
\half\int_0^1|E_{U\mu}[D_tf|\calF_t]|^2dt\right)\,,
$$
where $U\mu$ denotes the image of $\mu$ under $U$.
\end{theorem}
\nproof
Let $A_\la$, for $\la>0$ be defined as
$$
A_\la=\{\al\in  L_a^2(\mu,H):\,J(\al)\leq \la\}\,,
$$
by the $1$-convexity of $f$, $A_\la$ is a, non-empty, convex subset of
$L_a^2(\mu,H)$. Assume that $(\al_n,\,n\in \N)\subset A_\la$ converges
to some $\al$ in $ L_a^2(\mu,H)$, let $T_n=I_W+\al_n$ and
$T=I_W+\al$. From
Theorem \ref{thm-0}, 
$$
H(T_n\mu|\mu)\leq\half E[|\al_n|_H^2]\,,
$$
hence the sequence $(dT_n\mu/d\mu:\,n\in \N)$ is
uniformly integrable. From Lusin's lemma $(f\circ T_n,n\geq 1)$
converges in $L^0(\mu)$ to $f\circ T$. The
Fatou Lemma gives
$$
\al\geq\lim\inf_nE[f\circ T_n+\half|\al_n|_H^2]\geq E[f\circ T+\half|\al|_H^2]\,,
$$
i.e., $\al\in A_\la$, consequently $A_\la$ is closed in $ L_a^2(\mu,H)$, by convexity it is
weakly closed, which implies the weak lower semi continuity of
$J$. We claim that $A_\la$ is also bounded; in fact we have 
\beaa
\half\|\al\|_{L^2_a(\mu,H)}^2&=&J(\al)-E[f\circ T]\\
&\leq&J(\al)+E[f^-\circ T]\\
&\leq&J(\al)+E[e^{cf^-}]+\frac{1}{c}H(T\mu|\mu)\\
&\leq&J(\al)+E[e^{cf^-}]+\frac{1}{c}\|\al\|_{L^2_a(\mu,H)}^2\,,
\eeaa
where the last inequality follows again from Theorem \ref{thm-0} and
this estimate proves the boundedness of $A_\la$. Consequently $A_\la$  is a weakly
compact subset of $L^2_a(\mu,H)$, hence  $J$ attains its infimum on
it, convexity implies that this minimum is global. To see that $f$ is
a tamed functional under the last hypothesis, it suffices to remark
that this hypothesis and the assumption $E[e^{f^-}]<\infty$ imply that $E[(|f|+1)e^{-f}]<\infty$.
\nqed

\noindent
The hypothesis of integrability of $f$ seems to be indispensable for
the existence of a minimizing element $u$ as the following theorem
shows:
\begin{theorem}
\label{con-ex-thm}
Let $A\subset W$ be a measurable, $H$-convex set with positive Wiener
measure, let $f$ be defined as 
$$
e^{-f}=\mu(A) 1_A\,,
$$
usually $f$ is denoted by $\chi_A$ in analysis. 
Then there is no $\mu$-a.s. left inverible  API, say  $U=I_W+u$, 
which maps $W$ to $A$  unless $\mu(A)=1$. In other words $\chi_A$ is not a tamed Wiener
 functional for $\mu(A)\in (0,1)$.
\end{theorem}
\nproof
Suppose the contrary, then  there exists a $u\in L^2_a(\mu,H)$ such that $U=I_W+u$ is $\mu$-a.s. left
invertible and that $L=dU\mu/d\mu=e^{-f}/\mu(A)$, besides, from
Theorem \ref{attn-thm},  we should
have $J_\star=J(u)$.  Moreover, we can write 
$L=\rho(-\delta v)$, where $v\in L^2_a(U\mu,H)$, with $L\circ
U\,\rho(-\delta u)=1$ $\mu$-a.s. These equalities imply immediately
that 
$$
\rho(-\delta u)=\exp\left(-\delta u-\half|u|_H^2\right)=\mu(A)
$$
$\mu$-a.s. Consequently $\delta u=-\log \mu(A)-\half |u|_H^2$
$\mu$-a.s. Then the Kazamaki condition for the Girsanov exponential
(cf. \cite{BOOK}) of $u$ implies that 
\beaa
1&=&E\left[\exp\left(\la\delta
    u-\frac{\la^2}{2}|u|_H^2\right)\right]\\
&=&E\left[\exp\left(-\la\log\mu(A)-\frac{\la}{2}(1+\la)|u|_H^2\right)\right]
\eeaa
for any $\la\geq 0$ and this relation is possible only when $u=0$
$\mu$-a.s. and $\mu(A)=1$.
\nqed

\begin{remarkk}
Theorem \ref{con-ex-thm} gives an explicit  negative result about the
representability of positive random variables via API (cf.\cite{FUZ}),
although the same question has a positive answer if we forget about
the adaptedness, then there exists always an invertible perturbation of identity
$T=I_W+\nabla\phi$ which maps $W$ onto $A$, where $\phi\in \DD_{2,1}$ such that the operator
norm of  $\nabla^2\phi$ is essentially  bounded, cf.\cite{F-U3}.
\end{remarkk}
\begin{corollary}
\label{cor-1}
For any $H$-convex subset $A$,$\chi_A$ is tamed if and only if
$\mu(A)=1$
\end{corollary}

\noindent
The corollary below gives an explicit example of a stochastic
differential equation having a weak but no strong solution in a frame
which is totally different than the example of Tsirelson
(cf. \cite{I-W}):
\begin{corollary}
\label{cor-3}
Let $A\subset W$ be an $H$-convex set with $\mu(A)\in (0,1)$, define 
$$
d\nu=\frac{1}{\mu(A)}1_Ad\mu\,,
$$
and define $v$ as 
$$
\frac{1_A}{\mu(A)}=\rho(-\delta v)
$$
defined $\nu$-a.s., where $v\in L_a^0(\nu,H)$.  Let
$V_t=W_t+\int_0^t\dot{v}_sds$. Then under the probability $\nu$, $V$
is a Brownian motion and the canonical path $(W_t,t\in [0,1])$ is the weak
solution of the following stochastic differential equation
$$
W_t=-\int_0^t\dot{v}_s\circ W+V_t
$$
and this equation has no strong solution.
\end{corollary}
\nproof
The fact that $(W_t,t\in[0,1])$ is a weak solution follows from the fact that
$V$ is a Brownian motion under $\nu$, the fact that it is not a strong
solution follows from Theorem \ref{con-ex-thm}.\nqed

\section{\bf{Variational techniques to calculate the minimizers}}

\noindent
In this section we shall derive a necessary and sufficient condition
for a large class of adapted perturbation of identity. We begin with
some technical results:
\begin{lemma}
\label{var-1-lemma}
Assume that $g\in \DD_{p,1},\,p>1$ with $E[\exp\eps|\nabla
g|_H]<\infty$ for some $\eps>0$. Suppose that $\la\to\xi_\la$ is an
absolutely continuous curve from $[0,1]$ to $L^2_a(\mu,H)$ such that 
$$
\xi'_\la=\frac{d\xi_\la}{d\la}\in L^\infty_a(\mu,H). 
$$
Then we have
\begin{equation}
\label{basic-2}
g(w+\xi_\la(w))=g(w+\xi_0(w))+\int_0^\la(\nabla g(w+\xi_t(w)),\xi'_t(w))_Hdt
\end{equation}
$\mu$-almost surely.
\end{lemma}
\nproof 
Let $(P_t,t\geq 0)$ be the Ornstein-Uhlenbeck semigroup on $W$, let
$g_n=P_{1/n}g$, denote $w\to w+\xi_\la(w)$ by $T_\la(w)$. Since $g_n$
is an $H-C^\infty$-map (cf. \cite{BOOK}), (\ref{basic-2}) holds for
$g_n$. 
To pass to the limit, it suffices to show that $(|\nabla g_n\circ
T_\la|_H,n\geq 1)$ is uniformly integrable w.r.t. $d\la\times
d\mu=d\eta$. To see this, let $L_\la$ be the Radon-Nikodym density of
$T_\la\mu$ w.r.t. $\mu$ and  write
\beaa
E_\eta[|\nabla g_n\circ T_\cdot|_H 1_{\{|\nabla g_n\circ
  T_\cdot|_H>c\}}&=&E_\eta[|\nabla g_n|_H 1_{\{|\nabla g_n|_H>c\}}L_\cdot]\\
&\leq&E_\eta\left[e^{|\nabla g_n|}1_{\{|\nabla
    g_n|_H>c\}}\right]+E_\eta\left[\frac{1}{\eps}L_\cdot\log L_\cdot 1_{\{|\nabla
    g_n|_H>c\}}\right]\,.
\eeaa
The first term of the last line tends to zero uniformly in $n$ as
$c\to \infty$ due to the uniform integrability of $(\exp|\nabla
g_n|_H,n\geq 1)$, the second term also has the same behaviour by the
dominated convergence theorem, i.e., for any $\ga>0$, there exists
some $c_\ga>0$ such that 
$\sup_n\eta\{|\nabla g_n|_H>c\}<\beta$ as soon as $c>c_\ga$.
\nqed

\noindent
This lemma says that under its  regularity assumptions, to find the
candidate elements of $L^2_a(\mu,H)$ for the solution of the
minimization problems, we have to verify if they satisfy the
functional equation $u+\Phi(u)=0$, where $\Phi$ is 
defined by
$$
\Phi(\xi)=-\pi(\nabla f\circ (I_W+\xi))\,,
$$
and where $\pi$ denotes the dual predictable projection.
Suppose that $\|\nabla^2 f\|_{\rm op}\leq c<1$ almost surely,
where $c>0$ is a fixed constant and the norm is the operator norm
on $H$. Then the map $\Phi:L_{2}^a(\mu,H)\to L_{2}^a(\mu,H)$
is a strict
contraction, hence there exists a unique $u\in\DD_{2,0}^a(H)$
which satisfies the equation
$$
\dot{u}_t+E[D_tf\circ U|\calF_t]=0
$$
$dt\times d\mu$-almost surely. However, this condition is too
restrictive to be applicable and we reduce these hypothesis in the
next theorem.
Although the conditions are strong,  the
following theorem to justifies  rigourously the ideas explained in
the introduction and by itself it is an interesting property. Briefly
it tells that if  $f$ is sufficiently regular,
then every local minimum of $J$ is a global one, hence it is unique and the
corresponding API is almost surely invertible:
\begin{theorem}
\label{loc-glob-thm}
Assume that $f\in \DD_{p,1}$ for some $p>1$ and 
\begin{equation}
\label{hyp-f}
E[e^{\eps|f|}]<\infty\,.
\end{equation}
Let $L$ denote the
density $e^{-f}/E[e^{-f}]$. Define $v$ as to be 
$$
v(t,w)=\int_0^t\dot{v}_s(w)ds
$$
where 
$$
\dot{v}_t=E_L[D_tf|\calF_t]
$$
and where $E_L$ denotes the expectation operator w.r.t. the measure $Ld\mu$.
Assume that 
\begin{equation}
\label{hyp-v}
E[\exp\eps\|\nabla v\|^2_{op}+\exp\eps|v|_H^2]<\infty\,,
\end{equation}
for some arbitrary constant $\eps>0$, where $\|\cdot\|_{op}$ denotes the operator norm on
the Cameron-Martin space $H$. Let $u\in L^2_a(\mu,H)$ be a solution 
of the functional equation
$$
u+\Phi(u)=0\,.
$$
Then $u$ is a global minimizer of $J$, hence, in particular it is unique and
the conclusions of Theorem \ref{attn-thm} hold true for $U=I_W+u$.
\end{theorem}
\nproof
Let us note immediately that, using the Young inequality and  the
hypothesis (\ref{hyp-f}), (\ref{hyp-v}) imply  that $v\circ U\circ T$ is
again in $L^2_a(\mu,H)$ for any $T=I_W+\xi$ with $\xi\in
L^\infty_a(\mu,H)$ and the hypothesis of Lemma \ref{var-1-lemma} are
satisfied so that we have the relation 
$$
|v\circ (U+\la\xi)+\la \xi|_H^2=|v\circ U+u|_H^2+2\int_0^\la(\nabla
v\circ(U+t\xi)\xi+\xi,v\circ(U+t\xi)+t\xi)_Hdt\,, 
$$
which justifies the calculation of the directional derivative of $J$
which is done below.
v
\noindent
From the It\^o representation theorem, we have
$$
L=\exp\left(-\delta v-\half|v|_H^2\right)\,,
$$
hence 
$$
f=-\log E[e^{-f}]+\delta v+\half|v|_H^2\,.
$$
Using the last expression we get
\beaa
J(u)&=&E[f\circ U+\half|u|_H^2]\\
&=&-\log E[e^{-f}]+E[(\delta v+\half|v|_H^2)\circ U+\half|u|_H^2]\\
&=&-\log E[e^{-f}]+\half E[|v\circ U+u|_H^2]\,.
\eeaa
Note that the last equality above follows from the fact
$$
E[|v\circ U|_H^2]=E[|v|_H^2\,\frac{dU\mu}{d\mu}]\leq
E[e^{\eps|v|_H^2}]+\frac{1}{\eps}H(U\mu|\mu)<\infty\,,
$$
hence $E[\delta(v\circ U)]=0$.
Using Lemma \ref{var-1-lemma}, we obtain
\begin{equation}
\label{tan}
E[(\nabla v\circ U[\eta]+\eta,v\circ U+u)_H]=0
\end{equation}
for any $\eta\in L^2_a(\mu,H)$ bounded. Note that, for any $\xi\in L^2_a(\mu,H)$,
$$
 E[(\nabla v\circ U[\eta], \xi)_H]=E\int_0^1ds\int_s^1
 D_s\dot{v}_\tau\circ U \dot{\eta}_s\dot{\xi}_\tau d\tau\,.
$$
Since $\eta$ is adapted, $\tau\to \nabla_\eta\dot{v}_\tau$ is adapted,
consequently  $\nabla v\circ U$ is a quasi-nilpotent operator,
therefore $I_H+\nabla v\circ U$ is almost surely invertible. Relation (\ref{tan}) implies then  that $(I_H+\nabla
v^\star\circ U)(v\circ U+u)=0$ a.s, by the invertibility of
$I_H+\nabla v\circ U$, 
we obtain $v\circ U+u=0$ a.s., and this proves
the a.s. left invertibility of $U$ and the rest of the  proof follows from Theorem
\ref{attn-thm}.
\nqed


\vspace{2cm}

{\footnotesize{\bf{
\noindent
A.S. \"Ust\"unel, Telecom-Paristech (formerly ENST),
 Dept. Infres,\\
46, rue Barrault, 75013 Paris, France\\
email: ustunel@telecom-paristech.fr}
}}

\end{document}